\documentclass[reqno]{amsart}

\usepackage{amsmath,amssymb,amsthm, epsfig}
\usepackage{amsfonts}
\usepackage{graphicx}
\usepackage{amscd}
\usepackage{amsmath}
\usepackage{amssymb}
\usepackage[mathscr]{euscript}

\newtheorem{lemma}{Lemma}
\newtheorem{proposition}{Proposition}
\newtheorem{remark}{Remark}
\newtheorem{theorem}{Theorem}
\newtheorem{example}{Example}
\numberwithin{equation}{section}

\renewcommand{\div}{{\rm div}}
\newcommand{\R}{\mathring{R}}
\newcommand{\Rc}{\mathring{Ric}}

\makeatletter
\@namedef{subjclassname@2020}{%
  \textup{2020} Mathematics Subject Classification}
\makeatother

	\title[Compact quasi-Einstein manifolds with boundary]{Remarks on compact quasi-Einstein \\ manifolds with boundary}
	\author[R. Di\'ogenes]{R. Di\'ogenes}
	\author[T. Gadelha]{T. Gadelha}
	\author[E. Ribeiro Jr]{E. Ribeiro Jr}
	
	\address[R. Di\'ogenes]{UNILAB, Instituto de Ci\^encias Exatas e da Natureza, Rua Jos\'e Franco de Oliveira, 62790-970, Reden\c{c}\~ao - CE, Brazil.}\email{rafaeldiogenes@unilab.edu.br}
	
	\address[T. Gadelha]{Instituto Federal do Cear\'a - IFCE, Campus Maracana\'u, Av. Parque Central, 61939-140, Maracana\'u - CE, Brazil.}\email{tiago.gadelha@ifce.edu.br}

	\address[E. Ribeiro Jr]{Universidade Federal do Cear\'a - UFC, Departamento  de Matem\'atica, Campus do Pici, Av. Humberto Monte, 60455-760, Fortaleza - CE, Brazil.}
	\email{ernani@mat.ufc.br}
	
\thanks{T. Gadelha was partially supported by FUNCAP/Brazil}	
	
	 \thanks{E. Ribeiro was partially supported by CNPq/Brazil (\# 305410/2018-0 \& 160002/2019-2)}
	 
	 	\date{\today}
	
	

\begin{document}

	\begin{abstract} In this paper, we prove that a compact quasi-Einstein manifold $(M^n,\,g,\,u)$ of dimension $n\geq 4$ with boundary $\partial M,$ nonnegative sec\-tio\-nal curvature and  zero radial Weyl tensor is either isometric, up to scaling, to the standard hemisphere $\Bbb{S}^n_+,$ or $g=dt^{2}+\psi ^{2}(t)g_{L}$ and $u=u(t),$ where $g_{L}$ is Einstein with nonnegative Ricci curvature. A similar classification result is obtained by assuming a fourth-order vanishing condition on the Weyl tensor. Moreover, a new example is presented in order to justify our assumptions. In addition, the case of dimension $n=3$ is also discussed.
	\end{abstract}

	\maketitle
	
	\section{Introduction}
	\label{intro}

	 According to the works \cite{CaseShuWey, He-Petersen-Wylie2012,Petersen-Chenxu}, a complete $n$-dimensional Riemannian ma\-ni\-fold $(M^n,\,g),$ $n\geq 2,$ possibly with boundary $\partial M,$ is called an $m$-{\it quasi-Einstein manifold}, or simply {\it quasi-Einstein manifold}, if there exists a smooth potential function $u$ on $M^n$ obeying the following system
\begin{equation}
\label{eqdef}
\left\{%
\begin{array}{lll}
    \displaystyle \nabla^{2}u = \dfrac{u}{m}(Ric-\lambda g) & \hbox{in $M,$} \\
    \displaystyle u>0 & \hbox{on $int(M),$} \\
        \displaystyle u=0 & \hbox{on $\partial M,$} \\
    \end{array}%
\right.
\end{equation} for some constants $\lambda$ and $0<m<\infty.$ When $m=1$ we make the additional condition $\Delta u=-\lambda u.$ Here, $\nabla^{2} u$ stands for the Hessian of $u$ and $Ric$ is the Ricci tensor of $g.$

The study of $m$-quasi-Einstein manifolds is directly related to the existence of warped product Einstein metrics on a given manifold, see, e.g., \cite[pg. 265]{Besse} and \cite{Ernani2,Besse,CaseShuWey,MR,Rimoldi}. Another interesting motivation to investigate quasi-Einstein manifolds comes from the study of diffusion operators by Bakry and \'Emery \cite{bakry}, which is closely tied to the theory of smooth metric measure spaces (see also \cite{CaseT,MR,Rimoldi,WW}). Moreover, $1$-quasi-Einstein manifolds are more commonly called {\it static spaces}, which have attracted a large interest for their connections to the positive mass theorem and general relativity (see \cite[Remark 2.3]{CaseShuWey}). Hence, classifying $m$-quasi-Einstein manifolds or understanding their geometry is definitely an important problem.

Before proceeding, it is important to mention two examples of compact $m$-quasi-Einstein manifolds with boundary and constant scalar curvature (see \cite{DG 2019}).  The first one is the hemisphere $\Bbb{S}^n_+$ with the standard metric $g=dr^2+\sin^2r g_{\Bbb{S}^{n-1}}$ and  potential function $u(r)=\cos r,$ where $r$ is a height function with $r\leq\frac{\pi}{2}.$ The second one is $\Big[0,\sqrt{m/\lambda}\,\pi\Big]\times\Bbb{S}^{n-1},$ for $\lambda>0,$ with metric $g=dt^2+\frac{n-2}{\lambda}g_{\Bbb{S}^{n-1}}$ and potential function $u(t,x)=\sin\left(\sqrt{\lambda/m}\,t\right).$ Others examples of compact and noncompact quasi-Einstein manifolds can be found in, e.g., \cite{Besse,CaseP,CaseT,CaseShuWey,He-Petersen-Wylie2012,Ernani_Keti,Rimoldi,Wang}. In this article, we focus on nontrivial compact $m$-quasi-Einstein manifolds with boundary. Therefore, by the work \cite[Theorem 4.1]{He-Petersen-Wylie2012}, they have necessarily $\lambda>0.$

Inspired by some successful works on gradient Ricci solitons, locally conformally flat quasi-Einstein manifolds were investigated in \cite{CMMR,He-Petersen-Wylie2012}.  We refer to the reader to \cite{caoALM11} for an overview on Ricci solitons. Although many interesting results on gradient Ricci solitons are also obtained on quasi-Einstein manifolds, there exist examples of quasi-Einstein manifolds that are in stark contrast to the gradient Ricci solitons.  In this context, it is known that a locally conformally flat compact gradient Ricci soliton has constant curvature.  However, such a conclusion can not be extended to the compact quasi-Einstein manifolds (see \cite{Bo,He-Petersen-Wylie2012}). In \cite{He-Petersen-Wylie2012}, He, Petersen and Wylie showed that a compact quasi-Einstein manifold with harmonic Weyl tensor (i.e., $\div (W)=\nabla_{l}W_{ijkl}=0$) and $W(\nabla u, \cdot, \cdot, \nabla u)=0$ must have the form $$g=dt^2+\psi^{2}(t)g_{L}\,\,\,\,\,\hbox{and}\,\,\,\,\,u=u(t),$$ where $g_{L}$ is an Einstein metric and $(L,g_{L})$ has nonnegative Ricci curvature; see also \cite{Rio} and \cite{CatinoDGA} for the half-conformally flat case in dimension $4.$ As observed in \cite[Table 2]{He-Petersen-Wylie2012}, there exist examples of $(n\ge 5)$-dimensional quasi-Einstein manifolds which have harmonic Weyl tensor and $W(\nabla u, \cdot, \cdot, \nabla u)=0$ but are not locally conformally flat. Later,  Chen and He \cite{CC} proved that a Bach-flat  $(n\geq 4)$-dimensional compact quasi-Einstein manifold $M^n$ with boundary $\partial M$ is either Einstein or a finite quotient of a warped product with $(n-1)$-dimensional Einstein fiber.  We refer the readers to \cite{CatinoPJM} for a general discussion on Bach-flat Einstein-type manifolds.

To state our results, it is fundamental to remember that a Riemannian manifold $(M^n, \,g)$ has {\it zero radial Weyl curvature} when
\begin{equation}
\label{zeroWradial}
 W(\cdot, \,\cdot, \,\cdot,\nabla u) = 0,
 \end{equation} for a suitable potential function $u$ on $M^n. $ This class of manifolds clearly includes the case of locally conformally flat manifolds. In recent years, this condition has been combined with other geometric assumptions in order to obtain new classification results for gradient Ricci solitons, quasi-Einstein manifolds and critical metrics; see,  for instance, \cite{Baltazar,BRRemarks,CatinoZ, HPW2010,He-Petersen-Wylie2012}.

In \cite{BRRemarks}, Baltazar and Ribeiro showed that a compact $(n\geq 3)$-dimensional critical metric of the volume functional  $(M^n,\,g)$ with smooth boundary $\partial M,$ nonnegative sectional curvature and zero radial Weyl curvature is isometric to a geodesic ball in a simply connected space form $\Bbb{R}^n$ or $\Bbb{S}^n.$ These metrics are also known in the literature as $V$-{\it static spaces} or {\it Miao-Tam critical metrics} (see, e.g., \cite{BRRemarks,Batista,CEM,miaotam}).  In dimension $n=3,$ Huiya He \cite{HHe} improved this result obtained by Baltazar and Ribeiro by replacing the nonnegative sectional curvature by the nonnegative Ricci curvature condition.  Moreover, Huiya He obtained a similar result for $3$-dimensional static spaces (see also \cite{Lucas}). Therefore, since critical metrics of the volume functional and static spaces have necessarily constant scalar curvature, it is natural to ask whether a compact quasi-Einstein manifold with boundary, constant scalar curvature and nonnegative Ricci curvature must be Einstein.  However, we are able to exhibit an example of compact $(n\ge 4)$-dimensional quasi-Einstein manifold with boundary which has constant scalar curvature and nonnegative Ricci curvature, but it is not Einstein. To be precise, we have the following example.

\begin{example}
\label{exA}
	Consider $\Bbb{S}^{p+1}_+\times\Bbb{S}^q$ with the doubly warped product metric $$g=dr^2+\sin^2r g_{\Bbb{S}^p}+\frac{q-1}{p+m}g_{\Bbb{S}^q},$$ where $r(x,y)=h(x)$ and $h$ is a height function on $\Bbb{S}^{p+1}_+.$ Moreover, we consider  $u=\cos r$ with $r\leq\frac{\pi}{2}$ and $\lambda=p+m.$ 
	\end{example}
	
A detailed description of Example \ref{exA} will be presented in Section \ref{preliminaries}. This example raised the question of whether a compact quasi-Einstein manifold with boundary, constant scalar curvature and nonnegative sectional curvature is necessarily Einstein.  In our first result, we shall provide an answer to this question for dimension $n\geq 4$ by assuming the condition (\ref{zeroWradial}). More precisely, we have established the following result.

\begin{theorem}
	\label{thmbach-flat}
	Let $\big(M^{n},\,g,\,u,\,\lambda \big),$ $n\geq 4,$ be a nontrivial compact simply connected $(m\neq 1)$-quasi-Einstein manifold with smooth bo\-un\-dary $\partial M,$ constant scalar curvature and zero radial Weyl curvature.  Suppose that $M^n$ has nonnegative sec\-tio\-nal curvature. Then $\big(M^{n},\,g,\,u)$ is either 
	\begin{enumerate}
	\item isometric, up to scaling, to the standard hemisphere $\Bbb{S}^n_+,$ or
	\item $g=dt^{2}+\psi ^{2}(t)g_{L},$ $u=u(t),$ where $g_{L}$ is Einstein with nonnegative Ricci curvature.
	\end{enumerate}
\end{theorem}

\begin{remark}
Observe that Example \ref{exA} has nonnegative sectional curvature and constant scalar curvature, but it has no zero radial Weyl curvature. Therefore, the zero radial Weyl curvature condition (\ref{zeroWradial}) in Theorem \ref{thmbach-flat} can not be removed.
\end{remark}

\begin{remark}
It is known by \cite[Theorem 1.3]{Petersen-Chenxu} that every $3$-dimensional compact quasi-Einstein manifold $(M^{3},\,g,\,u)$ with boundary and constant scalar curvature is Einstein or its universal cover is a product of Einstein manifolds. 
\end{remark}

In another direction, Catino, Mastrolia and Monticelli \cite{catino} obtained an in\-te\-res\-ting classification result for gradient Ricci solitons admitting a {\it fourth-order vanishing condition on the Weyl tensor.} To be precise, they showed that any $(n \ge 4)$-dimensional gradient shrinking Ricci soliton with fourth-order divergence-free Weyl tensor (i.e., $\div^4(W)=\nabla_{j}\nabla_{k}\nabla_{l}\nabla_{i}W_{ijkl}=0$) is either Einstein or a finite quotient of $N^{n-k} \times \Bbb{R}^k,$ $(k > 0),$ i.e., the product of an Einstein manifold $N^{n-k}$ with the Gaussian shrinking soliton $\Bbb{R}^{k}.$ Notice that the fourth-order divergence-free Weyl tensor assumption is clearly weaker than locally conformally flat and harmonic Weyl tensor conditions consi\-de\-red in \cite{CMMR,He-Petersen-Wylie2012}. Moreover, this assumption was recently used  to classify critical metrics of the volume functional, static spaces and CPE metrics, for more details see, e.g., \cite{Baltazar,BDR,catino,Santos,YZ}. It is also important to emphasize that Example \ref{exA} has fourth-order divergence-free Weyl tensor, but it has no zero radial Weyl curvature.

 In the sequel, inspired by the work \cite{catino}, we establish a classification result for compact $(m\neq 1)$-quasi-Einstein manifold with smooth bo\-un\-dary under the fourth-order divergence-free Weyl tensor condition. More precisely, we have the following classification result.

	\begin{theorem}
	\label{div4W}
	Let $\big(M^{n},\,g,\,u,\,\lambda \big),$ $n\geq 4,$ be a nontrivial compact simply connected $(m\neq 1)$-quasi-Einstein manifold with smooth bo\-un\-dary $\partial M$ and zero radial Weyl curvature.  Suppose that $M^{n}$ has fourth-order divergence-free Weyl tensor $($i.e., $\div ^{4}(W)=0).$ Then $\big(M^{n},\,g,\,u)$ is either 
	\begin{enumerate}
	\item isometric, up to scaling, to the standard hemisphere $\Bbb{S}^n_+,$ or
	\item $g=dt^{2}+\psi ^{2}(t)g_{L},$ $u=u(t),$ where $g_{L}$ is Einstein with nonnegative Ricci curvature.
	\end{enumerate}
\end{theorem}

\begin{remark} Notice that Theorem \ref{div4W} improves \cite[Theorem 1.2]{He-Petersen-Wylie2012}  and \cite[Theorem 1.1]{CMMR} in the case that $M^n$ is compact with boundary.  Besides, Exam\-ple \ref{exA} guarantees that the zero radial Weyl curvature condition in Theorem \ref{div4W} can not be removed. 
\end{remark}

In dimension $n=3,$ we have established the following result.  

\begin{theorem}
\label{thmC}
Let $\big(M^{3},\,g,\,u,\,\lambda \big)$ be a nontrivial three-dimensional compact $m$-quasi-Einstein manifold with smooth bo\-un\-dary $\partial M.$ Suppose that $M^{n}$ has third-order divergence-free Cotton tensor $($i.e., $\div^{3}(C)=0). $ Then $(M^3,\,g)$ is locally conformally flat.
\end{theorem}

\begin{remark}
We point out that in dimension $n=3$ there are local solutions to the $m$-quasi-Einstein equations which are not locally conformally flat, see \cite[Example 3.5]{He-Petersen-Wylie2012}.
\end{remark}

	\section{Background}
	\label{preliminaries}

	In this section, we will present basic facts that will be useful for the establishment of the main results.  Moreover,  we will describe all details concerning Example \ref{exA} which motivated our results.  
	
	We start by recalling some special tensors in the study of curvature for a Riemannian manifold $(M^n,\,g)$ of dimension $n\ge 3.$  The first one is the Weyl tensor $W$ which is defined by the following decomposition formula
	\begin{eqnarray}\label{weyl tensor}
	R_{ijkl}&=&W_{ijkl}+\dfrac{1}{n-2}\big (R_{ik}g_{jl}+R_{jl}g_{ik}-R_{il}g_{jk}-R_{jk}g_{il}\big) \nonumber \\ 
	&&-\dfrac{R}{(n-1)(n-2)}\big (g_{jl}g_{ik}-g_{il}g_{jk}\big ),
	\end{eqnarray} where $R_{ijkl}$ denotes the Riemann curvature tensor. The second one is the Cotton tensor $C$ given by
	\begin{eqnarray}\label{cotton tensor}
	C_{ijk}=\nabla_{i}R_{jk}-\nabla_{j}R_{ik}-\dfrac{1}{2(n-1)}\big(\nabla_{i}Rg_{jk}-\nabla_{j}Rg_{ik}\big).
	\end{eqnarray} Notice that $C_{ijk}$ is skew-symmetric in the first two indices and trace-free in any two indices. The Weyl and Cotton tensors are related as follows
	\begin{equation}\label{cottonwyel}
	\displaystyle{C_{ijk}=-\frac{(n-2)}{(n-3)}\nabla_{l}W_{ijkl},}
	\end{equation} for $n\ge 4.$ In terms of the Schouten tensor 
	
	\begin{equation}
	\label{Schouten}
	A_{ij}=R_{ij}-\frac{R}{2(n-1)}g_{ij},
	\end{equation} we have
	
	\begin{equation}
	\label{Schouten2}
W_{ijkl}=R_{ijkl}-\frac{1}{(n-2)}\left(g_{ik}A_{jl}-g_{il}A_{jk}-g_{jk}A_{il}+g_{jl}A_{ik}\right)
	\end{equation} and
	
	\begin{equation}
	\label{Schouten3}
	C_{ijk}=\nabla_{i} A_{jk}-\nabla_{j}A_{ik}.
	\end{equation} In particular, it is well known that, for $n = 3,$ $W_{ijkl}$ vanishes identically, while $C_{ijk}=0$ if and only if $(M^3,\,g)$ is locally conformally flat; for $n\geq 4,$ $W_{ijkl}=0$ if and only if $(M^n,\,g)$ is locally conformally flat. 
	
The Bach tensor on a Riemannian manifold $(M^n,g),$ $n\geq 4,$ is defined in terms of the components of the Weyl tensor $W_{ikjl}$ as follows
	\begin{equation}
	\label{bach} B_{ij}=\frac{1}{n-3}\nabla_{k}\nabla_{l}W_{ikjl}+\frac{1}{n-2}R_{kl}W_{ikjl}.
	\end{equation} For $n=3,$ it is given by 
	\begin{equation}
	\label{Bach2}
B_{ij}=\nabla_{k}C_{kij}.
\end{equation} We say that $(M^n,\,g)$ is Bach-flat when $B_{ij}=0.$ In particular,  either locally conformally flat or Einstein metrics are necessarily Bach-flat. For $n\geq 4,$ it is known among the experts (see \cite[Lemma 5.1]{CaoChen}) the following formula for the divergence of the Bach tensor

	\begin{equation}
	\label{divbach}
	\nabla_{i}B_{ij}=\frac{n-4}{(n-2)^{2}}C_{jks}R_{ks}.
	\end{equation}

In order to proceed,  we remember that the fundamental equation of an $m$-quasi-Einstein manifold $(M^{n},\,g,\,u),$ possibly with boundary, is given by 
	\begin{equation}\label{fundamental equation}
	\nabla ^{2}u = \dfrac{u}{m}(Ric-\lambda g),
	\end{equation} where $u>0$ in the interior of $M^n$ and $u=0$ on $\partial M.$ Taking the trace of (\ref{fundamental equation}) one concludes that 
	\begin{equation}\label{laplaciano}
	\Delta u = \frac{u}{m}(R-\lambda n).
	\end{equation} Combining  Eqs. (\ref{fundamental equation}) and (\ref{laplaciano}) we get
	\begin{equation}\label{traceless Ricci}
	u\, \mathring{Ric} = m \mathring{\nabla^2}\, u,
	\end{equation} where $\mathring{T}=T-\dfrac{{\rm tr}T}{n}g$ stands for the traceless part of $T.$ Besides, we recall that the scalar curvature of a compact quasi-Einstein manifold with boundary must satisfy

	\begin{equation}
	\label{eq3e4}
	R\geq\frac{n(n-1)}{m+n-1}\lambda.
	\end{equation} For more details, see \cite[Remark 5.1]{He-Petersen-Wylie2012} and \cite[Proposition 3.6]{CaseShuWey}.
	
	For any $m$-quasi-Eisntein manifold, it turns out that the Cotton tensor $C_{ijk}$ can be expressed in terms of the Weyl tensor $W_{ijkl}$ and an auxiliary tensor $T_{ijk}$ (for the proof see \cite{CC,Ranieri_Ernani}).

\begin{lemma}[\cite{CC,Ranieri_Ernani}]
\label{CWT}
Let $\big(M^{n},\,g\big)$ be a Riemannian manifold satisfying (\ref{eqdef}) for some potential function $u$ on $M^n.$ Then, it holds
\begin{equation}
uC_{ijk}=mW_{ijkl}\nabla_lu+T_{ijk},
\end{equation}
where $T$ is given by 
\begin{eqnarray}\label{T-tensor}
T_{ijk}&=&\frac{m+n-2}{n-2}(R_{ik}\nabla_ju-R_{jk}\nabla_iu)+\frac{m}{n-2}(R_{jl}\nabla_lug_{ik}-R_{il}\nabla_lug_{jk})\nonumber\\
 &&+\frac{(n-1)(n-2)\lambda+mR}{(n-1)(n-2)}(\nabla_iug_{jk}-\nabla_jug_{ik})\nonumber\\&&-\frac{u}{2(n-1)}(\nabla_iRg_{jk}-\nabla_jRg_{ik}).
\end{eqnarray}
\end{lemma} 

Moreover, as a direct consequence of (\ref{T-tensor}) one obtains that
\begin{eqnarray}\label{LemmaCijk.Tijk}
C_{ijk}\nabla_{i}uR_{jk}&=&-\frac{n-2}{2(m+n-2)}C_{ijk}T_{ijk}.
\end{eqnarray} For more details, see Eq. (2.13) in \cite{DG 2019}.

We further recall a B\"ochner type formula for $m$-quasi-Einstein manifolds established in \cite[Theorem 3]{DG 2019}.

\begin{proposition}[\cite{DG 2019}]
\label{th3}
Let $\big(M^{n},\,g,\,u,\,\lambda \big)$ be an $m$-quasi-Einstein manifold. Then we have
\begin{eqnarray*}
\frac{1}{2}\div \big(u\nabla|\mathring{Ric}|^{2}\big)&=&u|\nabla\mathring{Ric}|^{2}+\frac{m(n-2)u}{m+n-2}|C|^2+u\langle\nabla^2R,\Rc\rangle\nonumber\\&&+\frac{m+2n-2}{n-1}\Rc(\nabla u,\nabla R)+\frac{m-1}{2}\langle\nabla|\Rc|^2,\nabla u\rangle\nonumber\\
	&&+\frac{2Ru}{n-1}|\Rc|^2+\frac{2nu}{n-2}{\rm tr}(\Rc^3)-2uW_{ikjp}\R_{ij}\R_{kp}\nonumber\\
	&&-\frac{m^2(n-2)}{m+n-2}C_{ijk}W_{ijkl}\nabla_lu\nonumber,
\end{eqnarray*} where ${\rm tr}(\Rc^3)=\R_{ij}\R_{jk}\R_{ki}$ and $\mathring{Ric}=Ric-\frac{R}{n}g.$
\end{proposition}

In the sequel we shall describe all details concerning Example \ref{exA}. Before to do so, we recall a basic result that can be found in \cite[Section 2.4, pg. 71]{petersen}. For a Riemannian manifold $I\times\Bbb{S}^p\times\Bbb{S}^q$ endowed with metric $$g=dr^2+\varphi^2(r)g_{\Bbb{S}^p}+\psi^2(r)g_{\Bbb{S}^q},$$ it holds
	\begin{eqnarray}\label{Hessr}
	\nabla^2r=\varphi'\varphi g_{\Bbb{S}^p}+\psi'\psi g_{\Bbb{S}^q}.
	\end{eqnarray} Moreover, for any $X\in\mathfrak{X}(\Bbb{S}^p)$ and $Y\in\mathfrak{X}(\Bbb{S}^q),$ we have
	\begin{eqnarray}\label{RicCampos}
	Ric(\partial_r)&=&\left(-p\frac{\varphi''}{\varphi}-q\frac{\psi''}{\psi}\right)\partial_r,\nonumber\\
	Ric(X)&=&\left(-\frac{\varphi''}{\varphi}+(p-1)\frac{1-\varphi'^2}{\varphi^2}-q\frac{\varphi'\cdot\psi'}{\varphi\cdot\psi}\right)X,\\
	Ric(Y)&=&\left(-\frac{\psi''}{\psi}+(q-1)\frac{1-\psi'^2}{\psi^2}-p\frac{\varphi'\cdot\psi'}{\varphi\cdot\psi}\right)Y.\nonumber
	\end{eqnarray}

	Now, we are ready to discuss Example \ref{exA}.
	
	\begin{example}[Example \ref{exA}]
		Consider $M=\Bbb{S}^{p+1}_+\times\Bbb{S}^q$ with the doubly warped product metric $$g=dr^2+\sin^2rg_{\Bbb{S}^p}+\frac{q-1}{p+m}g_{\Bbb{S}^q},$$ where $r(x,y)=h(x)$ and $h$ is the height function on $\Bbb{S}^{p+1}_+.$ In particular, $r\leq\frac{\pi}{2}.$ So, by using (\ref{Hessr}) we have $\nabla^2r=\sin r\cos rg_{\Bbb{S}^p}.$ Next, we set $u=\cos r$ and $\lambda=p+m.$ Hence, $\nabla u=-\sin r\nabla r$ and
		\begin{eqnarray*}
			\nabla^2u&=&-\sin r\nabla^2r-\cos rdr^2\\
			&=&-\sin^2r\cos rg_{\Bbb{S}^p}-\cos rdr^2\\
			&=&-u(dr^2+\sin^2rg_{\Bbb{S}^p}).
		\end{eqnarray*} Thereby, for $X\in\mathfrak{X}(\Bbb{S}^{p})$ and $Y\in\mathfrak{X}(\Bbb{S}^q)$ we deduce
		\begin{eqnarray}\label{Ex4-Hess}
		\nabla^2u(\partial_r,\partial_r)&=&-udr^2(\partial_r,\partial_r)=-ug(\partial_r,\partial_r),\nonumber\\
		\nabla^2u(X,X)&=&-u\sin^2rg_{\Bbb{S}^p}(X,X)=-ug(X,X),\\
		\nabla^2u(Y,Y)&=&0.\nonumber
		\end{eqnarray}
		
		On the other hand, it follows from (\ref{RicCampos}) that
		\begin{eqnarray}
		Ric(\partial_r)&=&-p\frac{(-\sin r)}{\sin r}\partial_r=p\partial_r,\nonumber\\
		Ric(X)&=&\left(1+(p-1)\frac{1-\cos^2r}{\sin^2r}\right)X=pX,\\
		Ric(Y)&=&(q-1)\frac{p+m}{q-1}Y=(p+m)Y.\nonumber
		\end{eqnarray} Therefore, considering $E=\frac{u}{m}(Ric-\lambda g)=\frac{u}{m}(Ric-(p+m)g),$ we infer
		\begin{eqnarray}\label{Ex4-T}
		E(\partial_r,\partial_r)&=&\frac{u}{m}\left[Ric(\partial_r,\partial_r)-(p+m)g(\partial_r,\partial_r)\right]\nonumber\\
		&=&-ug(\partial_r,\partial_r),\nonumber\\
		E(X,X)&=&\frac{u}{m}\left[Ric(X,X)-(p+m)g(X,X)\right]\\
		&=&-ug(X,X),\nonumber\\
		E(Y,Y)&=&\frac{u}{m}\left[Ric(Y,Y)-(p+m)g(Y,Y)\right]=0.\nonumber
		\end{eqnarray} Finally, combining (\ref{Ex4-Hess}) and (\ref{Ex4-T}), we then obtain $$\nabla^2u=\frac{u}{m}(Ric-\lambda g).$$ Moreover, $u=\cos\frac{\pi}{2}=0$ on $\partial M.$ This shows that $(M,\,g,\,u,\lambda)$ is a nontrivial compact $m$-quasi-Einstein manifold with boundary.
	\end{example}

	\section{Proof of the Main Results}
	\label{secA}
	
	In this section we are going to present the proofs of Theorems \ref{thmbach-flat}, \ref{div4W} and \ref{thmC}.

	\subsection{The proof of Theorem \ref{thmbach-flat}}

	\begin{proof} To begin with, taking into account that $M^{n}$ has constant scalar curvature and zero radial Weyl tensor, we may invoke Proposition \ref{th3} to get the following divergence formula
\begin{eqnarray}\label{corWeylZeroeq1}
\frac{1}{2}\div \big(u\nabla|\mathring{Ric}|^{2}\big)&=&u|\nabla\mathring{Ric}|^{2}+\frac{m(n-2)u}{m+n-2}|C|^2+\frac{m-1}{2}\langle\nabla|\Rc|^2,\nabla u\rangle\nonumber\\
&&+\frac{2Ru}{n-1}|\Rc|^2+\frac{2nu}{n-2}{\rm tr}(\Rc^3)-2uW_{ikjp}\R_{ij}\R_{kp},
\end{eqnarray}  where ${\rm tr}(\Rc^3)=\mathring{R_{ij}}\mathring{R_{jk}}\mathring{R_{ki}}.$

On the other hand,  the decomposition of the Riemann tensor (\ref{weyl tensor}) on a Riemannian manifold $(M^n,\,g)$ gives

\begin{eqnarray*}
	R_{ij}R_{jk}R_{ik}-R_{ijkl}R_{jl}R_{ik}&=&\frac{n}{n-2}R_{ij}R_{jk}R_{ik}-W_{ijkl}R_{jl}R_{ik}\nonumber\\&&-\frac{(2n-1)}{(n-1)(n-2)}R|Ric|^{2}+\frac{R^{3}}{(n-1)(n-2)},
\end{eqnarray*} so that

\begin{eqnarray}
\label{12fg}
R_{ij}R_{jk}R_{ik}-R_{ijkl}R_{jl}R_{ik}&=&\frac{n}{n-2}R_{ij}R_{jk}R_{ik}-W_{ijkl}R_{jl}R_{ik}\nonumber\\&&-\frac{(2n-1)}{(n-1)(n-2)}R|\mathring{Ric}|^{2}-\frac{1}{n(n-2)}R^{3}.\nonumber\\
\end{eqnarray} Besides, one easily verifies that

$$R_{ij}R_{ik}R_{jk}=\mathring{R}_{ij}\mathring{R}_{jk}\mathring{R}_{ik}+\frac{3}{n}R|\mathring{Ric}|^{2}+\frac{R^{3}}{n^{2}}$$ and therefore, (\ref{12fg}) becomes

\begin{equation}\label{eq2thmbach-flat}
\frac{n}{n-2}tr(\mathring{Ric}^{3})-W_{ijkl}R_{ik}R_{jl}=R_{ij}R_{jk}R_{ik}-R_{ijkl}R_{jl}R_{ik}-\frac{1}{n-1}R|\mathring{Ric}|^{2}. 
\end{equation}

At the same time, since $M^n$ has constant scalar curvature, we can use \cite[Proposition 3.3]{Petersen-Chenxu} (see also \cite[Lemma 3.2]{CaseShuWey}) to deduce
\begin{eqnarray*}
|\Rc|^2=-\frac{m+n-1}{n(m-1)}(R-n\lambda)\left(R-\frac{n(n-1)}{m+n-1}\lambda\right).
\end{eqnarray*} Hence, $|\Rc|^2$ is also constant.  Now, substituting (\ref{eq2thmbach-flat}) into (\ref{corWeylZeroeq1}) and using that $|\Rc|^2$ is constant we arrive at

\begin{eqnarray}\label{eq4thmbach-flat}
0&=&u|\nabla\mathring{Ric}|^{2}+\frac{m(n-2)u}{m+n-2}|C|^2+2u\left(R_{ij}R_{jk}R_{ik}-R_{ijkl}R_{jl}R_{ik}\right).
\end{eqnarray} In order to proceed, we need to guarantee that the last term of the right hand side of (\ref{eq4thmbach-flat}) is nonnegative.  Indeed, choosing $\lambda_{i}$ to be the eingenvalues of the Ricci curvature, we obtain 
		
$$\left(\nabla_{i}\nabla_{j}R_{ik}-\nabla_{j}\nabla_{i}R_{ik} \right) R_{jk}=\sum_{i<j} R_{ijij}(\lambda_{i}-\lambda_{j})^{2}.$$ Then, taking into account that $M^n$ has nonnegative sectional curvature we deduce 

\begin{equation}
\label{pkl1}
\left(\nabla_{i}\nabla_{j}R_{ik}-\nabla_{j}\nabla_{i}R_{ik}\right)R_{jk}\geq 0.
\end{equation}

Proceeding,  from the commutation formulae (Ricci identities) for any Riemannian manifold $(M^{n},\,g),$ it holds
\begin{equation}\label{idRicci}
\nabla_i\nabla_j R_{kl}-\nabla_j\nabla_i R_{kl}=R_{ijks}R_{sl}+R_{ijls}R_{ks},
\end{equation}  which combined with (\ref{pkl1}) yields 

$$R_{ij}R_{jk}R_{ik}-R_{ijkl}R_{jl}R_{ik}\geq 0.$$ This data into (\ref{eq4thmbach-flat}) allows us to infer the following inequality

		\begin{eqnarray}\label{eq5thmbach-flat}
		0\geq u|\nabla \Rc|^{2} + \frac{m(n-2)}{m+n-2} u|C|^{2}.
		\end{eqnarray} Whence,  since $g$ and $u$ are real analytic in harmonic coordinates (see, e.g.,  \cite[Proposition 2.4]{He-Petersen-Wylie2012}) and taking into account that $u$ is nonnegative, it follows from (\ref{eq5thmbach-flat}) that the Cotton tensor $C$ vanishes and $\nabla \Rc =0.$

		To conclude,  we use (\ref{bach}), (\ref{fundamental equation}) and (\ref{cottonwyel})  in order to obtain

\begin{eqnarray*}
(n-2)uB_{ij}&=& uW_{ikjl}R_{kl}\nonumber\\&=& mW_{ikjl}\nabla_{k}\nabla_{l}u\nonumber\\&=&m\nabla_{k}(W_{ijkl}\nabla_{l}u)-m\nabla_{k}W_{ikjl}\nabla_{l}u\nonumber\\&=&-m\nabla_{k}W_{jlik}\nabla_{l}u \nonumber\\&=&m\left(\frac{n-3}{n-2}\right)C_{jli}\nabla_{l}u,
\end{eqnarray*}  so that

$$u B_{ij} = m\frac{(n-3)}{(n-2)^2}C_{jli}\nabla_{l}u=0,$$ where we have used that the Cotton tensor $C$ vanishes. This implies that $(M^{n},g)$ is Bach-flat and hence,  it suffices to invoke \cite[Corollary 1.7]{CC} by Chen and He to conclude that $M^{n}$ is either Einstein or $$g=dt^{2}+\psi ^{2}(t)g_{L},\,\,\,\,\, \,\,\,\,\,u=u(t),$$ where $g_{L}$ is Einstein with nonnegative Ricci curvature. In the first case, we may apply  Proposition 2.4 of \cite{Petersen-Chenxu} to conclude that $M^n$ is isometric, up to scaling, to the standard hemisphere $\Bbb{S}^n_+.$  So, the proof is completed.
	\end{proof}

		\subsection{The proof of Theorem \ref{div4W}}
		
		\begin{proof}
		Initially, we invoke (\ref{divbach}) to infer
		
	\begin{eqnarray}
	(n-2)\int_{M}u\, \div^{2}(B) dV_{g}&=&\frac{(n-4)}{(n-2)}\int_{M}u \nabla_{i}(C_{ijk}R_{jk}) dV_{g} \nonumber \\
	&=&\frac{(n-4)}{(n-2)}\int_{M}\nabla_{i}\left(u C_{ijk}R_{jk}\right)dV_{g}\nonumber\\&&-\frac{(n-4)}{(n-2)}\int_{M} \nabla_{i}uC_{ijk}R_{jk} dV_{g}, \nonumber
		\end{eqnarray} and by using (\ref{LemmaCijk.Tijk}) one sees that
		
			\begin{eqnarray}\label{eq1proofdiv3C}
(n-2)\int_{M}u\, \div^{2}(B) dV_{g}&=&\frac{(n-4)}{2(m+n-2)}\int_{M} C_{ijk}T_{ijk} dV_{g}. 
	\end{eqnarray}

On the other hand, it follows from (\ref{cottonwyel}) and (\ref{bach}) that $$(n-2)B_{ij}=\nabla_{k}C_{kij}+W_{ikjl}R_{kl}$$ and therefore, our assumption that $M^n$ has zero radial Weyl curvature implies that

\begin{eqnarray}\label{eq2proofdiv3C}
		(n-2)\int_{M}u\, \div^{2}(B) dV_{g} &=&\int_{M}u\, \nabla_{j}\nabla_{i}\nabla_{k}C_{kij} dV_{g} + \int_{M}u\nabla_{j}\nabla_{i}\left(W_{ikjl}R_{kl}\right) dV_{g} \nonumber \\
		&=&\int_{M}u \div^{3}(C)  dV_{g} + \int_{M}\nabla_{j}\left(u \nabla_{i}(W_{ikjl}R_{kl})\right) dV_{g} \nonumber \\
		&&-\int_{M} \nabla_{j}u\nabla_{i}(W_{ikjl}R_{kl})  dV_{g} \nonumber \\ &=&\int_{M}u\,\div^{3}(C)  dV_{g}-\int_{M}\nabla_{j}u\nabla_{i}W_{ikjl}R_{kl} dV_{g}. 
\end{eqnarray} Now, substituting  (\ref{cottonwyel}) into (\ref{eq2proofdiv3C}) yields

	\begin{eqnarray}\label{eq3proofdiv3C}
		(n-2)\int_{M}u\, \div^{2}(B) dV_{g}&=&\int_{M}u\, \div^{3}(C) dV_{g}-\frac{(n-3)}{n-2}\int_{M} C_{jlk}\nabla_{j}uR_{kl} dV_{g}. \nonumber \\
			\end{eqnarray} and by using  (\ref{LemmaCijk.Tijk}), one obtains that
			
			\begin{eqnarray}\label{eq4proofdiv3C}
		(n-2)\int_{M}u\, \div^{2}(B) dV_{g}&=&\int_{M}u\,  \div^{3}(C) dV_{g}+\frac{(n-3)}{2(m+n-2)}\int_{M} C_{jlk}T_{jlk} dV_{g}. \nonumber \\
		\end{eqnarray} This jointly with (\ref{eq1proofdiv3C}) gives

			\begin{equation*}
			\int_{M}u\,\div^{3}(C) dV_{g}+\frac{1}{2(m+n-2)}\int_{M} T_{ijk}C_{ijk} dV_{g}=0.			
			\end{equation*} Next,  using once more that $M^n$ has zero radial Weyl curvature and Lemma \ref{CWT} we get   
				\begin{eqnarray*}
	\int_{M}u\, \div^{3}(C) dV_{g}+\frac{1}{2(m+n-2)}\int_{M}u |C|^{2}dV_{g}=0.
	\end{eqnarray*} Now, taking into account that  $M^{n}$ has fourth-order divergence-free weyl tensor,  which implies from (\ref{cottonwyel}) that $\div^{3}(C)=0,$ one sees that
	\begin{eqnarray*}
	\int_{M}u |C|^{2}dV_{g}=0.
	\end{eqnarray*} Therefore, the Cotton tensor $C$ vanishes.  From now on, it suffices to apply the same arguments as in final steps of the proof of Theorem \ref{thmbach-flat} to infer that $(M^{n},g)$ is Bach-flat and hence,  we invoke Corollary 1.7 of \cite{CC} combined with Proposition 2.4 of \cite{Petersen-Chenxu} to conclude that $(M^{n},g)$ is either isometric, up to scaling, to the standard hemisphere $\Bbb{S}^n_+,$ or $g=dt^{2}+\psi ^{2}(t)g_{L},$ $u=u(t),$ where $g_{L}$ is Einstein with nonnegative Ricci curvature. This finishes the proof of the theorem.						
	\end{proof}

\subsection{Proof of Theorem \ref{thmC}}

\begin{proof}
We start by claiming that 

\begin{equation}
\label{eqBC1}
\nabla_{j}B_{ij}=-R_{jk}C_{ijk},
\end{equation} (see also Eq. (4.2) in \cite{Caosteady}). Indeed,  by combining (\ref{Schouten3}) and (\ref{Bach2}) one sees that 
$$B_{ij}=\nabla_{k}\left(\nabla_{k}A_{ij}-\nabla_{i}A_{kj}\right)$$ and consequently, 

\begin{eqnarray}
\nabla_{i}B_{ij}&=&\nabla_{i}\nabla_{k}\left(\nabla_{k}A_{ij}-\nabla_{i}A_{kj}\right)\nonumber\\&=&\left(\nabla_{i}\nabla_{k}-\nabla_{k}\nabla_{i}\right)\nabla_{k}A_{ij}\nonumber\\&=& -R_{il}\nabla_{l}A_{ij}+R_{kl}\nabla_{k}A_{lj}+R_{ikjl}\nabla_{k}A_{il}\nonumber\\&=& R_{ikjl}\nabla_{k}A_{il}.
\end{eqnarray} Then, it suffices to use (\ref{Schouten2}) and (\ref{Schouten3}) to infer 
\begin{eqnarray*}
\nabla_{i}B_{ij}&=&\left(g_{ij}A_{kl}-g_{il}A_{kj}-g_{kj}A_{il}+g_{kl}A_{ij}\right)\nabla_{k}A_{il}\nonumber\\&=& A_{jk}g_{il}C_{lki}+A_{ik}C_{kji}=-R_{ki}C_{jki},
\end{eqnarray*} which proves (\ref{eqBC1}). 

Proceeding,  one easily verifies from (\ref{eqBC1}) that

\begin{eqnarray*}
u\nabla_{j}\nabla_{i}B_{ij}&=&-u\nabla_{j}\left(R_{ki}C_{jki}\right)\nonumber\\&=& -\nabla_{j}\left(uR_{ki}C_{jki}\right)+\left(\nabla_{j}u\right)R_{ki}C_{jki},\nonumber
\end{eqnarray*}  and then by using (\ref{LemmaCijk.Tijk}) and Lemma \ref{CWT}, one obtains that

\begin{eqnarray}
u\nabla_{j}\nabla_{i}B_{ij}&=&-\nabla_{j}\left(uR_{ki}C_{jki}\right)-\frac{1}{2(m+1)}C_{jki}T_{jki}\nonumber\\&=&-\nabla_{j}\left(u R_{ki}C_{jki}\right)-\frac{1}{2(m+1)}u|C|^{2}.
\end{eqnarray} Upon integrating this expression over $M^3,$ we use Stokes' theorem and our assumption that $\div^{3} (C)=\div^{2} (B)=0$ to infer $$\int_{M} u \div^{2} (B) dV_{g}=-\int_{M}u|C|^{2} dV_{g}=0$$ and hence,  $C$ vanishes identically. So, the proof is completed.
\end{proof}

\end{document}